\documentclass{amsart}

\usepackage[all]{xy}

\usepackage{pictexwd, epsfig, amsthm, amssymb, amsmath, graphicx, psfrag,amsxtra, latexsym}
\usepackage{enumerate}

\newcommand{\mZ}{\mathbb Z}
\newcommand{\mR}{\mathbb R}
\newcommand{\tK}{\tilde K}

 \DeclareMathOperator{\z}{{\mathbb
Z}} 
\DeclareMathOperator{\g}{\Gamma}

 \DeclareMathOperator{\coker}{coker}

\DeclareMathOperator{\f}{\mathcal F}
\DeclareMathOperator{\V}{\mathcal V}
\DeclareMathOperator{\tf}{\tilde {\mathcal F}}
\DeclareMathOperator{\kz}{\mathbb K \mathbb Z^{-\infty}}
\DeclareMathOperator{\kr}{\mathbb K R^{-\infty}}
\DeclareMathOperator{\vc}{\mathcal V\mathcal C}

\DeclareMathOperator{\fin}{\mathcal F\mathcal I\mathcal N}
\DeclareMathOperator{\all}{\mathcal A\mathcal L\mathcal L}

\theoremstyle{plain}
\newtheorem{Thm}{Theorem}

\newtheorem{Cor}[Thm]{Corrollary}

\newtheorem*{main}{Main Theorem}

\theoremstyle{definition}

\theoremstyle{remark}
\newtheorem*{Prf}{Proof}

\begin{document}

\title[Splitting formulas for certain Waldhausen Nil-groups.]
      {Splitting formulas for certain Waldhausen Nil-groups.}
\author{Jean-Fran\c{c}ois\ Lafont}
\address{Department of Mathematics\\
         Ohio State University\\
         Columbus, OH  43210}
\email[Jean-Fran\c{c}ois\ Lafont]{jlafont@math.ohio-state.edu}
\author{Ivonne J.\ Ortiz}
\address{Department of Mathematics and Statistics\\
         Miami University\\
         Oxford, OH 45056}
\email[Ivonne J.\ Ortiz]{ortizi@muohio.edu}

\begin{abstract}
We provide splitting formulas for certain Waldhausen Nil-groups.  We
focus on Waldhausen Nil-groups associated to {\it acylindrical} amalgamations
$\g=G_1*_HG_2$ of groups $G_1,G_2$ over a common subgroup $H$.  
For these amalgamations, we explain how, provided $G_1,G_2,\g$ satisfy the
Farrell-Jones isomorphism conjecture, the Waldhausen Nil-groups
$Nil^W_*(R H; R[G_1-H], R[G_2-H])$ can be expressed as a direct sum
of Nil-groups associated to a specific collection of virtually
cyclic subgroups of $\g$. A special case covered by our theorem is
the case of arbitrary amalgamations over a finite group $H$.  
\end{abstract}

\maketitle

\section{Introduction}

Waldhausen's Nil-groups were introduced in the two foundational
papers \cite{W1}, \cite{W2}. The motivation behind these Nil-groups
originated in a desire to have a Mayer-Vietoris type sequence in
algebraic K-theory.  More precisely, if a group $\g=G_1*_HG_2$
splits as an amalgamation of two groups $G_1,G_2$ over a common
subgroup $H$, one can ask how the algebraic $K$-theory of the group
ring $R \g$ is related to the algebraic K-theory of the integral
group rings $R G_1$, $R G_2$, $R H$.  Motivated by the corresponding
question in homology (or cohomology), one might expect a
Mayer-Vietoris type exact sequence:
$$\ldots \rightarrow K_{i+1}(R \g) \rightarrow K_{i}(R H)\rightarrow K_{i}(R G_1)\oplus
K_i(R G_2) \rightarrow K_{i}(R \g) \rightarrow \ldots $$ A major
result in \cite{W1}, \cite{W2} was the realization that the
Mayer-Vietoris sequence above holds, provided one inserts suitable
``error-terms'', which are precisely the Waldhausen Nil-group
associated to the amalgamation $\g=G_1*_HG_2$.  In general,
associated to any ring $S$ (such as $R H$), and any pair of flat
$R$-bimodules $M_1,M_2$ (such as the $R [G_i-H]$), Waldhausen
defines Nil-groups $Nil^W_*(S; M_1,M_2)$. The Waldhausen Nil-groups
$Nil^W_*(R H; R[G_1-H], R[G_2-H])$ are precisely the ``error-terms''
mentioned above.

\vskip 5pt

Another context in which these Nil-groups make an appearance has to
do with the {\it reduction to finites}.  To explain this we recall
the existence of a generalized equivariant homology theory, having
the property that for any group $\Gamma$, one has an isomorphism:
$$H_n^\Gamma(*; \kr) \cong K_n(R \Gamma).$$
The term appearing to the right is the homology of the
$\Gamma$-space consisting of a point $*$ with the trivial
$\Gamma$-action.  Now for any $\Gamma$-space $X$, the obvious map
$X\rightarrow *$ is clearly $\Gamma$-equivariant, and hence induces
a homomorphism:
$$H_n^\Gamma (X; \kr) \rightarrow H_n^\Gamma(*; \kr) \cong K_n(R \Gamma).$$
The Farrell-Jones isomorphism conjecture \cite{FJ93} asserts that if
$X= E_{\vc}\Gamma$ is a model for the classifying space for
$\Gamma$-actions with isotropy in the family of virtually cyclic
subgroups, then the homomorphism described above is actually an
isomorphism.  Explicit models for $E_{\vc}\Gamma$ are known for only
a few classes of groups: virtually cyclic groups, crystallographic
groups (Alvez and Ontaneda \cite{AO06} and Connolly), hyperbolic
groups (Juan-Pineda and Leary \cite{JL06} and L\"uck \cite{Lu05}),
and relatively hyperbolic groups (Lafont and Ortiz \cite{LO}).  In
contrast, classifying spaces for {\it proper} actions, denoted
$E_{\fin}\Gamma$ are known for many classes of groups.  Now for any
group $\Gamma$, one always has a unique (up to $\Gamma$-equivariant
homotopy) map $E_{\fin}\Gamma \rightarrow E_{\vc}\Gamma$, which
induces a well-defined relative assembly map:
$$H_n^{\g} (E_{\fin}\Gamma; \kr) \rightarrow H_n^{\g} (E_{\vc}\Gamma;
\kr).$$ In view of the Farrell-Jones isomorphism conjecture, it is
reasonable to ask whether this latter map is itself an isomorphism.

In many cases, this map is known to be split injective.  For
instance in the case where $\g$ is a $\delta$-hyperbolic group (see
Rosenthal-Sch\"utz \cite{RS05}), or more generally, when $\g$ has
finite asymptotic dimension (see Bartels-Rosenthal \cite{BR}), the
relative assembly map is a split injection for arbitrary rings $R$.
From the topological viewpoint, the case where $R=\mZ$ is the most
interesting.  In this situation, Bartels \cite{Ba03} has shown that
the relative assembly map is split injective for arbitrary groups
$\g$.  If the relative assembly map discussed above is actually an
isomorphism, then we say that $\g$ satisfies the {\it reduction to
finites}.

Let us now specialize to the case where $R=\mZ$, i.e. we will be
focusing on integral group rings.  In this situation, the
obstruction to the relative assembly map being an isomorphism lies
in the Nil groups associated to the various infinite virtually
cyclic subgroups of $\Gamma$. More precisely, the following two
statements are equivalent (see \cite[Thm. A.10]{FJ93}):
\begin{enumerate}
\item The group $\g$ satisfies the reduction to finites.
\item Every infinite virtually cyclic subgroup $V \leq \g$ satisfies
the reduction to finites.
\end{enumerate}
Now for a virtually cyclic group $V$, the {\it failure} of the
reduction to finites can be measured by the cokernel of the relative
assembly map.  Let us recall that infinite virtually cyclic groups
$V$ come in two flavors:
\begin{itemize}
\item groups that surject onto the infinite dihedral group
$D_\infty$, and hence can be decomposed $V=A*_CB$, with $A,B,C$
finite, and $C$ of index two in the groups $A,B$. \item groups that
do {\it not} surject onto $D_\infty$, which can always be written in
the form $V=F \rtimes _\alpha \mZ$, where $F$ is a finite group and
$\alpha\in Aut(F)$.
\end{itemize}
In the case where $V$ surjects onto $D_\infty$, the cokernel of the
relative assembly map coincides with the Waldhausen Nil-group
associated to the splitting $V=A*_CB$.  In the case where $V$ does
not surject onto $D_\infty$, the cokernel of the relative assembly
map consists of two copies of the Farrell Nil-group associated to
$V=F \rtimes _\alpha \mZ$, denoted $NK_*(\mZ F, \alpha)$.

We now have two contexts in which Waldhausen Nil-groups make an
appearance: (1) they measure failure of the Mayer-Vietoris sequence
in algebraic $K$-theory, and (2) they contain obstructions for
groups to satisfy the reduction to finites.  Having motivated our
interest in these groups, we can now state our:

\begin{main}
Let $\g= G_1*_HG_2$ be an {\em acylindrical} amalgamation, and assume that the
Farrell-Jones Isomorphism Conjecture holds for the groups $\g,
G_1,G_2$. Denote by $\V$ a
collection of subgroups of $\g$ consisting of one representative $V$
from each conjugacy class of subgroups satisfying:
\begin{enumerate}
\item $V$ is virtually cyclic,
\item $V$ is {\it not} conjugate to a subgroup of $G_1$ or $G_2$, and
\item $V$ is maximal with respect to subgroups satisfying (1) and (2).
\end{enumerate}
Then for arbitrary rings $R$ we have an isomorphism for $*\leq 1$:
$$Nil^W_*(R H; R[G_1-H], R[G_2-H])\cong \bigoplus _{V\in \V}
H_*^{V}(E_{\fin}V\rightarrow *; \kr)$$ where
$H_*^{V}(E_{\fin}V\rightarrow *;\kr)$ are the cokernels of the
relative assembly maps associated to the virtually cyclic subgroups
$V\in \V$ (and hence consist of classical Farrell or Waldhausen
Nils).
\end{main}

\noindent{\bf Remarks: (1)} The notion of an {\em acyclindrical amalgamation}
was formulated by Sela \cite{Sel97} in relation to his work on the accessibility
problem for finitely generated groups.  An amalgamation
$\g= G_1*_HG_2$ is said to be acylindrical if there exists an
integer $k$ such that, for every path $\eta$ of length $k$ in the Bass-Serre
tree $T$ associated to the splitting of $\g$, the stabilizer of $\eta$
is finite.  Observe that if the amalgamating subgroup $H$ is finite,
the amalgamation is automatically acylindrical (with $k=1$).

\vskip 5pt

\noindent {\bf (2)}  The cokernels $H_*^{V_i}(E_{\fin}V_i\rightarrow
*;\kr)$ are the familiar Waldhausen or (two copies of the) Farrell Nil
groups, according to whether the virtually cyclic group $V_i$
surjects onto $D_\infty$ or not. Note that every virtually cyclic
subgroup $V$ that maps onto $D_\infty$ contains a canonical index
two subgroup $V^\prime$ which does {\it not} map onto $D_\infty$
(the pre-image of the obvious index two $\mZ$ subgroup in
$D_\infty$).  Recent independent work by various authors (Davis
\cite{D}, Davis-Khan-Ranicki \cite{DKR}, Quinn, and Reich),
has established that the Waldhausen Nil-group of $V$ is isomorphic
to the Farrell Nil-group of $V^\prime$.

\vskip 5pt

\noindent {\bf (3)} From the computational viewpoint, the Main Theorem
combined with the previous remark
completely reduces (modulo the Isomorphism Conjecture) the
computation of Waldhausen Nil-groups associated to acylindrical
amalgamations to that of Farrell Nil-groups.

\vskip 5pt

\noindent {\bf (4)} Consider the simple case 
of a {\it free product} $\g = G_1* G_2$.  In this
situation, the group $\g$ is known to be (strongly) relatively hyperbolic, relative
to the subgroups $G_1,G_2$.  Assuming the Farrell-Jones isomorphism conjecture 
for $\g$, previous work of the authors \cite[Cor. 3.3]{LO2} yields the following 
expression for $K_n(R \g)$:
$$H_n^{\g}(E_{\fin}\g) \oplus \Big( \bigoplus_{i=1,2} H_n^{G_i}(E_{\fin}G_i
\rightarrow E_{\vc}G_i)\Big) \oplus \Big(\bigoplus _{V\in \V} H_n^{V}(E_{\fin}V
\rightarrow *)\Big)$$
where $\V$ is the collection of virtually cyclic subgroups mentioned in our Main
Theorem (we omitted the coefficients $\kr$ to simplify notation).  Now, morally 
speaking, the Waldhausen Nil-group is the portion of the $K$-theory of $\g$ 
that does {\it not come from the $K$-theory of the factors $G_i$}.  Recalling
the well known fact that every finite subgroup of $\g$ has to be 
conjugate into one of the $G_i$, one sees that in the expression above, this
should just be the last term.  Our Main Theorem came about from trying
to make this heuristic precise.

\vskip 10pt

\section{Proof of Main Theorem}

Given the group $\g=G_1*_HG_2$ satisfying the hypotheses of our
theorem, let us form the family  $\f$ of subgroups of $\g$
consisting of all virtually cyclic subgroups that can be conjugated
into either $G_1\leq \g$ or $G_2\leq \g$.  Observe that we have a
containment of families $\f \subset \vc$ of subgroups of $\g$, which
in turn induce an assembly map:
$$\rho: H_n^\Gamma (E_{\f}\g; \kr) \rightarrow H_n^\Gamma (E_{\vc}\g; \kr).$$
Our proof will focus on analyzing the map $\rho$, and in particular
on gaining an understanding of the cokernel of that map.  Let us
start by describing the Waldhausen Nil-group as the cokernel of a
suitable assembly map. Corresponding to the splitting $\g =
G_1*_HG_2$, we have a simplicial action of $\g$ on the corresponding
Bass-Serre tree $T$ (see \cite{Ser80}). From the natural
$\g$-equivariant map $T\rightarrow *$, we get an assembly map:
$$\rho^\prime: H_n^\Gamma (T; \kr) \rightarrow H_n^\Gamma (*; \kr) \cong K_*(R \g).$$
We start out with the important:

\vskip 5pt

\noindent {\bf Fact:}  The map $\rho ^\prime$ is split injective,
and
$$\coker(\rho^\prime) \cong Nil^W_*(R H; R[G_1-H], R[G_2-H]).$$

\vskip 5pt

A proof of this {\bf Fact} can be found in Davis \cite[Lemma 7]{D} (see 
also Remark (i) at the end of this section).
In view of this result, we are merely trying to identify the
cokernel of the map $\rho^\prime$.  The first step is to relate the
cokernel of $\rho^\prime$ with the cokernel of the map $\rho$.

\vskip 5pt

\noindent {\bf Claim 1:}  The map $\rho$ is split injective, and there
is a canonical isomorphism
$$\coker(\rho)\cong \coker(\rho ^\prime)$$

\vskip 5pt

We observe that we have four families of subgroups that we are
dealing with: the three we have looked at so far are $\vc$, $\all$,
and the family $\f$ we introduced at the beginning of our proof
(consisting of virtually cyclic subgroups conjugate into one of the
$G_i$). In addition, there is the family $\mathcal{G}$ consisting of
{\it all} subgroups of $\g$ which can be conjugated into either
$G_1$ or $G_2$.  Now observe that we have containments of families
$\f \subset \mathcal{G} \subset \all$, and $\f \subset \vc \subset
\all$. Furthermore, we have that $T$ is a model for $E_{\mathcal
G}\g$. This yields the following commutative diagram:
$$\xymatrix{
H_*^\Gamma (E_{\f}\g ; \kr) \ar[r]^{\rho} \ar[d] &  H_*^\Gamma
(E_{\vc}\g ; \kr) \ar[d]^{\cong}\\
H_*^\Gamma (T ; \kr) \ar[r]_{\rho^\prime} &  H_*^\Gamma (* ; \kr) \\
}$$ where all the maps are relative assembly maps corresponding to
the inclusions of the various families of subgroups.  Note that the
horizontal maps are precisely the ones we are trying to relate.  Now
recall that we are assuming that $\g$ satisfies the Farrell-Jones
isomorphism conjecture.  This immediately implies that the second
vertical map is an isomorphism, as indicated in the commutative
diagram.  So in order to identify the cokernels of the two
horizontal maps, we are left with showing that the first vertical
map is also an isomorphism.

The first vertical map is a relative assembly map, corresponding to
the inclusion of the families $\f \subset \mathcal{G}$ of subgroups
of $\g$.  In order to show that the relative assembly map is an
isomorphism, one merely needs to establish that for every maximal
subgroup $H\in \mathcal {G} - \f$, the corresponding relative
assembly map induced by the inclusions of families $\f (H) \subset
\mathcal{G}(H)$ of subgroups of $H$ is an isomorphism (see
\cite[Thm. A.10]{FJ93}, \cite[Thm. 2.3]{LS}). But observe that the maximal subgroups
in $\mathcal {G} - \f$ are precisely the (conjugates of) the
subgroups $G_i\leq \g$. Furthermore, for these subgroups, we have
that $\mathcal{G}(G_i)=\all (G_i)$, and that $\f (G_i) = \vc (G_i)$.
Hence the relative assembly maps which we require to be isomorphisms
are exactly those induced by $E_{\vc}G_i\rightarrow E_{\all}G_i
\cong *$, i.e. those that arise in the Farrell-Jones isomorphism
conjecture.  Since we are assuming the isomorphism conjecture holds
for the groups $G_1,G_2$, we conclude that the first vertical map is
indeed an isomorphism, completing the proof of the claim.

\vskip 10pt

At this point, combining {\bf Claim 1} with the {\bf Fact}, we have an
identification:
$$\coker(\rho)\cong \coker(\rho ^\prime)\cong  Nil^W_*(R H; R[G_1-H], R[G_2-H]).$$
In order to complete the proof, we now focus entirely on studying
the map $\rho$, with a goal of showing that one can express its
cokernel as a direct sum of the desired Nil-groups associated to the virtually
cyclic subgroups $V\in \V$. We
remind the reader that $\rho$ is the relative assembly map
induced by the map $E_{\f}\g \rightarrow E_{\vc}\g$,
where $\f$ is the family of subgroups consisting of all virtually
cyclic subgroups of $\g$ that can be conjugated into either $G_i$.

In order to analyze this relative assembly map, we will need
to make use of some properties of the $\g$-action on the
Bass-Serre tree.  Particularly, we would like to understand
the behavior of virtually cyclic subgroups $V\in \vc -\f$.
The specific result we will require is contained in our:

\vskip 10pt

\noindent {\bf Claim 2:} In the case of an acylindrical
amalgamation, the stabilizer of any geodesic $\gamma$ in the
Bass-Serre tree $T$ is a virtually cyclic subgroup of $\g$.
Furthermore, every virtually cyclic subgroup $V\leq \g$ satisfying
$V\in \vc - \f$ stabilizes a unique geodesic $\gamma \subset T$.

\vskip 10pt

Let us start by recalling some basic facts concerning
the action of $\g$ on the Bass-Serre tree $T$ corresponding to the
amalgamation $\g = G_1*_HG_2$:
\begin{itemize}
\item the action is without inversions, i.e. if an element stabilizes
an edge $e$, then it automatically preserves the chosen orientation
of $e$,
\item the stabilizer of any vertex $v\in T$ is isomorphic to a
conjugate of $G_1$ or $G_2$,
\item the stabilizer of any edge $e\subset T$ is isomorphic to a
conjugate of $H$,
\item any finite subgroup of $\g$ fixes a vertex in
$T$,
\end{itemize}
The first three statements above are built into the definition of
the Bass-Serre tree (see \cite{Ser80}), while the last statement 
is a well known general facts about group actions on trees.  
We remind the reader that
a {\em geodesic} in a tree $T$ will be a subcomplex simplicially
isomorphic to $\mR$, with the standard simplicial structure (i.e.
vertices at the integers, and edges between).

To show the first statement in our  claim, we note that
$Stab_{\g}(\gamma)$
clearly fits into a short exact sequence:
$$0\rightarrow Fix_{\g}(\gamma) \rightarrow Stab_{\g}(\gamma) \rightarrow
Sim_{\g , \gamma}(\mR) \rightarrow 0$$ where $Fix_{\g}(\gamma)$ is
the subgroup fixing $\gamma$ pointwise, and $Sim_{\g , \gamma}(\mR)$
is the induced simplicial action on $\mR$ (obtained by simplicially
identifying $\gamma$ with $\mR$).  Note that the group of simplicial
automorphisms of $\mathbb R$ is $D_\infty$, the infinite dihedral
group. In particular, we see that $Sim_{\g , \gamma}(\mR)$ is
virtually cyclic (in fact is isomorphic to $\mZ/2, \mZ$, or
$D_\infty$).

Next we observe that $Fix_{\g}(\gamma)$ is finite.  To see this, we
recall that the amalgamation $\g= G_1*_HG_2$ was assumed to be {\em
acylindrical}, which means that there exists an integer $k\geq 1$
with the property that the stabilizer $Stab_{\g}(\eta) \leq \g$ of
any combinatorial path $\eta\subset T$ of length $\geq k$ is finite.
Since $\gamma \subset T$ is a geodesic, it contains combinatorial
subpaths $\eta$ of arbitrarily long length (in particular, length
$\geq k$). The obvious containement $Fix_{\g}(\gamma) \leq
Stab_{\g}(\eta)$ now completes the argument for the first statement
in our claim.

For the second statement, we note that $\g$ acts simplicially on the
Bass-Serre tree, and hence the given virtually cyclic subgroup $V\in
\vc - \f$ likewise inherits an action on $T$. Since $V \notin \f$,
we have that the $V$ action on $T$ has no globally fixed point.  In
particular, $V$ must be an infinite virtually cyclic subgroup, and
hence contains elements of infinite order.  If $g\in V$ is an
arbitrary element of infinite order, we now claim that $g$ cannot
fix any vertex in $T$.

Assume, by way of contradiction, that there exists a vertex $v$
fixed by $g$ (and hence by $V^\prime$).  Let $T^\prime \subset T$ be
the subset consisting of points that are fixed by $V^\prime$.  Note
that $T^\prime$ is non-empty (since $v\in T^\prime$), and is a
subtree of $T$ (since $\g$ acts simplicially on $T$).  Furthermore,
observe that the group $F:=V/V^\prime$ inherits a simplicial action
on $T^\prime$.  But note that $F$ is finite, and hence the
$F$-action on $T^\prime$ has a fixed vertex $w\in T^\prime \subset
T$. But this immediately implies that the original group $V$ fixes
$w$, a contradiction as $V\notin F$.

Now establishing that $V$ stabilizes a geodesic is a straightforward
application of standard techniques in the geometry of group actions
on trees (applied to $T$). For the convenience of the reader, we
give a quick outline of the argument. For an arbitrary element $g$
of infinite order in $V$, one can look at the associated
displacement function on $T$, i.e. the distance from $v$ to $g\cdot
v$.  The previous paragraph establishes that this function is
strictly positive.  One then considers the set $Min(g)$ of points in
$T$ which minimize the displacement function, call this minimal
value $\mu _g$. It is easy to see that:
\begin{itemize}
\item $Min(g)$ contains a geodesic $\gamma$: take a vertex $v\in
Min(g)$, and consider $\gamma := \bigcup _{i\in \mZ}g^i\cdot \eta$,
where $\eta$ is the geodesic segment from $v$ to $g\cdot v$ (note
that such a non-trivial segment exists by the previous paragraph),
\item in fact, $\gamma =Min(g)$: any point
at distance $r>0$ from $\gamma$ will be displaced $2r + \mu_g>
\mu_g$, and so cannot lie in $Min(g)$,
\item for any non-zero integer $i$, $Min(g)=Min(g^i)$: any point
at distance $r>0$ from $\gamma$ will be displaced $2r+|i|\cdot \mu_g
> |i|\cdot \mu_g$, while points on $\gamma$ will clearly only be
displaced $|i|\cdot \mu_g$ by the element $g^i$,
\item for any two elements $g,h$ of infinite order in $V$, we have
$Min(g)=Min(h)$: two such elements have a common power, and apply
the previous statement.
\end{itemize}
From the observations above, we see that every single element in $V$
of infinite order stabilizes the exact same geodesic $\gamma \subset
T$.

Hence the only elements we might have to worry about are elements
$h\in V$ of {\it finite order}.  For these, we just note that $V$
contains $V^\prime \triangleleft V$, a finite index cyclic normal
subgroup generated by an element $g$ of infinite order.  We have a
natural morphism from $H=\langle h\rangle$ to $Aut(V^\prime)\cong
\mZ/2$.  In particular, we have that $hgh^{-1}=g^{\pm 1}$ and hence
we have the obvious equalities $d(g\cdot hv, hv)=d(h^{-1}gh \cdot v,
v)=d(g^{\pm 1} \cdot v,v)=\mu_{g^{\pm 1}}=\mu_g$. Since $hv$ is
minimally displaced by $g$, it must lie on $\gamma=Min(g)$.  This
deals with elements of finite order, hence completes the
verification that $\gamma$ is $V$-invariant.  Finally, from the 
fact that $V$ has a finite index subgroup that acts on $\gamma$
via a translation, it is easy to see that there are no other $V$-invariant
geodesics in $T$, yielding uniqueness.  This finishes the
argument for our {\bf Claim 2}.

\vskip 10pt

Having established some basic properties of the $\g$-action on $T$,
we now return to the main argument.  Recall that by combining
our {\bf Claim 1} with the {\bf Fact}, we have reduced the proof of
the Main Theorem to understanding the cokernel of the relative
assembly map:
$$\rho: H_n^\Gamma (E_{\f}\g; \kr) \rightarrow H_n^\Gamma (E_{\vc}\g; \kr)$$
where $\f$ is the family of subgroups consisting of all virtually
cyclic subgroups of $\g$ that can be conjugated into either $G_i$.

Now recall that in \cite{LO}, the authors introduced the notion of a
collection of subgroups to be {\it adapted} to a nested family of
subgroups, and in the presence of an adapted family, showed how a
model for the classifying space with isotropy in the smaller family
could be ``promoted'' to a model for the classifying space with
isotropy in the larger family.  In a subsequent paper \cite{LO2},
the authors used some recent work of Lueck-Weiermann \cite{LW} to
give an alternate model for this classifying space, which had the
additional advantage of providing {\it explicit splittings} for the
cokernel of the relative assembly maps.  Let us briefly recall the
relevant definitions.

Given a nested pair of families $\f \subset \tf$ of subgroups of
$\g$, we say that a collection $\{H_{\alpha}\}_{\alpha \in I}$ of
subgroups of $\g$ is {\it adapted} to the pair $(\f, \tf)$ provided
that:
\begin{enumerate}
\item For all $H_1, H_2 \in \{H_{\alpha}\}_{\alpha \in I}$,
either $H_1=H_2$, or $H_1 \cap H_2 \in \f$.
\item The collection $\{H_{\alpha}\}_{\alpha \in I}$ is {\it conjugacy
closed} i.e.\ if $H \in \{H_{\alpha}\}_{\alpha \in I}$ then
$gHg^{-1} \in \{H_{\alpha}\}_{\alpha \in I}$ for all $g \in \g$.
\item Every $H \in \{H_{\alpha}\}_{\alpha \in I}$ is {\it self-normalizing}, i.e.\ $N_{\Gamma}(H)=H$.
\item For all $G \in \tf \setminus \f$, there
exists $H \in \{H_{\alpha}\}_{\alpha \in I}$ such that $G \leq H$.
\end{enumerate}
In \cite{LO2}, we applied this result to the nested family $\fin
\subset \vc$ for relatively hyperbolic groups (for which an adapted
collection of subgroups is easy to find).

\vskip 5pt

In our present context, we would like to find a collection of
subgroups adapted to the nested pair of families $\f \subset \vc$.
Recall that $\f$ consists of all virtually cyclic subgroups that can
be conjugated into either $G_i$, and $\vc$ consists of all virtually
cyclic subgroups of $\g$.

\vskip 10pt

\noindent {\bf Claim 3:}  The collection $\{H_\alpha\}$ of subgroups
of $\g$ consisting of all {\it maximal} virtually cyclic subgroups
in $\vc - \f$ is adapted to the pair $(\f,\vc)$.

\vskip 10pt

To verify this, we first note that properties (2) and (4) in the
definition of an adapted family are immediate.  Property (3)
follows easily from {\bf Claim 2}: let $V\in
\{H_\alpha\}$ be given, and consider $N_{\g}(V)$.  We know that $V$
leaves invariant a unique geodesic $\gamma \subset T$.  Furthermore, for
every $g\in \g$, we see that $gVg^{-1}$ leaves $g\cdot \gamma$ invariant.
Uniqueness of the $V$-invariant geodesic $\gamma$ now implies that
$\gamma$ is actually $N_{\g}(V)$-invariant.  In particular, we have
containments $V\leq N_{\g}(V) \leq Stab_{\g}(\gamma)$.  But from
{\bf Claim 2}, we know that $Stab_{\g}(\gamma)$ is virtually cyclic, and
maximality of $V$ now forces all the containments to be equalities,
and in particular, $V= N_{\g}(V)$ as required by property (3).

For property (1), let $V_1,V_2 \in \{H_\alpha\}$.  We want to
establish that either $V_1=V_2$, or that $V_1\cap V_2 \in \f$.  So
let us assume that $V_1\neq V_2$.  We know from {\bf Claim 2} that 
each $V_i$ stabilizes a unique geodesic $\gamma_i$, and from the
maximality of the groups $V_i$, we actually have
$V_i=Stab_{\g}(\gamma_i)$.  Since $V_1\neq V_2$, we have that
$\gamma_1\neq \gamma_2$.  There are now two possibilities: (i)
either $\gamma_1\cap \gamma_2 = \emptyset$, or (ii) $\gamma_1 \cap
\gamma_2$ is a path in $T$.  We claim that in both cases, the intersection
$H=V_1\cap V_2 \leq V_i$ has the property that the $H$-action on the
corresponding $\gamma_i$ {\it fixes a point}.

To see this, let us first consider possibility (i): since $\gamma_1\cap
\gamma_2 =\emptyset$, one can consider the (unique) minimal length
geodesic segment $\eta$ joining $\gamma_1$ to $\gamma_2$.  We observe
that, since $H$ stabilizes both $\gamma_i$, it must leave the segment
$\eta$ invariant.  In particular, $H$ must fix the vertex $v_i=\eta \cap
\gamma_i \in \gamma_i$, as desired.
Next consider possibility (ii): if $\gamma_1\cap \gamma_2 \neq \emptyset$,
then the intersection will be a subpath $\eta \subset \gamma_i$.  Note that
$\eta$ is either a geodesic {\it segment}, or is a geodesic {\it ray}, and
in both cases, will be invariant under the group $H$.  It $\eta$
is a geodesic ray (i.e. homeomorphic to $[0,\infty)$), then there is a
(topologically) distinguished point inside $\eta$, which will have to
be fixed by $H$.  If $\eta$ is a geodesic segment, then each element
in $H$ either fixes $\eta$, or reverses $\eta$ (note that the latter can only
occur if $\eta$ has {\it even} length, as $\g$ acts on $T$ without
inversions).  In particular, we see that if $\eta$ has odd length, then
every point in $\eta$ is fixed by $H$, while if $\eta$ has even length,
then the (combinatorial) midpoint is fixed.

Finally, we observe that $H\leq V_i$ acts on $\gamma_i$, and fixes
a point.  This immediately implies that $H$ contains a subgroup of
index at most two which {\it acts trivially} on $\gamma_i$, i.e. $H^\prime
\leq Fix_{\g}(\gamma_i)$.  But recall that the latter group is finite
(see the proof of {\bf Claim 2}), completing the proof of property (1).
We conclude that the collection $\{H_\alpha\}$ is an adapted collection
of subgroups for the nested families $(\f, \vc)$, as desired.

\vskip 10pt

Finally, we exploit the adapted family we just constructed to establish:

\vskip 10pt

\noindent {\bf Claim 4:}  We have an identification:
$$\coker(\rho)\cong \bigoplus _{V\in \V} H_*^V(E_{\fin}V\rightarrow *; \kr)$$
where the groups $H_*^V(E_{\fin}V\rightarrow *; \kr)$ are the
cokernels of the relative assembly maps associated with the
virtually cyclic groups $V\in \V$

\vskip 10pt

The argument for this is virtually identical to the one given in 
\cite[Cor. 3.2]{LO2}; we reproduce the argument here for the 
convenience of the reader.
From {\bf Claim 3}, we have an adapted family for the pair $(\f , \vc)$,
consisting of all maximal subgroups in $\vc - \f$.  From this
adapted collection of subgroups, \cite[Prop. 3.1]{LO2} establishes
(using \cite[Thm. 2.3]{LW}) a method for constructing an
$E_{\vc}\g$. Namely, if $\V$ is a complete set of
representatives of the conjugacy classes within the adapted
collection of subgroups $\{H_\alpha\}$, form the cellular
$\g$-pushout:
$$\xymatrix{\coprod_{V \in \V} \g \times_V E_{\f}(V)   \ar[d] _\alpha \ar[rrr]^\beta & & &
E_{\f}(\g) \ar[d]\\
\coprod_{V\in \V} \g \times_V E_{\vc}(V) \ar[rrr] & & & X}
$$
Then the resulting space $X$ is a model for $E_{\vc}(\Gamma)$ (we
refer the reader to \cite[Prop. 3.1]{LO2} for a more precise
discussion of this result, including a description of the maps $\alpha, \beta$ in
the above cellular $\g$-pushout).  Note that the map $\rho$ whose
cokernel we are trying to understand is precisely the map on
(equivariant) homology induced by the second vertical arrow in the
above cellular $\g$-pushout.

Since $X$ is the double mapping cylinder of the maps $\alpha,\beta$
in the above diagram, one has a natural $\g$-equivariant
decomposition of $X$ by taking $A$ (respectively $B$) to be the
$[0,2/3)$ (respectively $(1/3,1]$) portions of the double mapping
cylinder.  Applying the homology functor $H_*^{\g}(-;\kr)$ (and
omitting the coefficients to shorten notation), we have the
Mayer-Vietoris sequence:
$$\ldots \rightarrow H_*^{\g}(A\cap B) \rightarrow H_*^{\g}(A)\oplus H_*^{\g}(B)
\rightarrow H_*^{\g}(X)\rightarrow H_{*-1}^{\g}(A\cap B)\rightarrow
\ldots$$ But now observe that we have obvious $\g$-equivariant
homotopy equivalences between:
\begin{enumerate}
\item $A \simeq _{\g} \coprod_{V\in \V} \g \times_V
E_{\vc}(V)$
\item $B\simeq _{\g} E_{\f}(\g)$
\item $A\cap B \simeq _{\g} \coprod_{V\in \V} \g \times_V
E_{\f}(V)$
\end{enumerate}
Furthermore, the homology theory we have takes disjoint unions into
direct sums.  Combining this with the induction structure, we obtain
isomorphisms:
$$H_*^{\g}(A) \cong \bigoplus _{V\in \V} H_*^{\g}(\g
\times_V E_{\vc}(V)) \cong \bigoplus _{V\in \V}
H_*^{V}(E_{\vc}(V))$$
$$H_*^{\g}(A\cap B) \cong \bigoplus _{V\in \V} H_*^{\g}(\g
\times_V E_{\f}(V)) \cong \bigoplus _{V\in \V}
H_*^{V}(E_{\f}(V))$$ Now observing that the groups $V\in \V$
are all virtually cyclic, we have that each $E_{\vc}(V)$ can be
taken to be a point.  Furthermore, for the groups $V\in \V$,
we have that the restriction of the family $\f$ to $V$ coincides
with the family of finite subgroups of $V$, i.e. $E_{\f}(V)=
E_{\fin}(V)$.  Substituting all of this in the above Mayer-Vietoris
sequence, we get the long exact sequence:
$$\ldots \rightarrow \bigoplus _{V\in \V}
H_*^{V}(E_{\fin}(V)) \rightarrow H_*^{\g}(E_{\f}\g)\oplus \bigoplus
_{V\in \V} H_*^{V}(*)\rightarrow H_*^{\g}(E_{\vc}\g)\rightarrow
\ldots$$ Now observe that the each of the maps
$H_*^{V}(E_{\fin}(V))\rightarrow H_*^{V}(*)$ are split injective
(see for instance \cite{RS05}). Since the map $\rho:
H_*^{\g}(E_{\f}\g) \rightarrow H_*^{\g}(E_{\vc}\g)$ is also split
injective (from {\bf Claim 1}), we now have an identification:
$$\coker (\rho) \cong \bigoplus _{V\in \V}
\coker \big( H_*^{V}(E_{\fin}(V)) \rightarrow H_*^{V}(*)\big)$$

\vskip 10pt

Finally, combining the {\bf Fact} with {\bf Claim 1} and {\bf Claim 4},
we see that we have the desired splitting:
$$Nil^W_*(R H; R[G_1-H], R[G_2-H])\cong \bigoplus _{V\in \V} H_*^V(E_{\fin}V\rightarrow *;\kr)$$
where the groups $H_*^{V}(E_{\fin}(V) \rightarrow *;\kr)$ denote the
cokernels appearing in {\bf Claim 4}. This completes the proof of
the Main Theorem.

\vskip 10pt

\noindent {\bf Remark: (1)} One of the key ingredients in our proof was
the {\bf Fact} established by Davis in \cite[Lemma 7]{D}.  Prior to
learning of Davis' preprint, the authors had an alternate argument
for the {\bf Fact}.  For the sake of the interested readers, we
briefly outline our alternate approach.

Anderson-Munkholm \cite[Section 7]{AM00} defined a functor,
continuously controlled $K$-theory, from the category of diagrams of
holink type to the category of spectra.  Munkholm-Prassidis
\cite[Theorem 2.1]{MP01} showed that the Waldhausen Nil-group we are
interested in can be identified with the cokernel of a natural split
injective map $\tK _{*+1}^{cc}(\xi ^+)\rightarrow K_*(\z \g)$, where
$\xi ^+$ is a suitably defined diagram of holink type associated to
the splitting $\g= G_1*_HG_2$ (see \cite[Section 9]{AM00}).
Furthermore, Anderson-Munkholm have shown \cite[Theorem 9.1]{AM00}
that there is a natural isomorphism $\tK _{*+1}^{cc}(\xi ^+) \cong
\tK _{*+1}^{bc}(\xi ^+)$, where the latter is the boundedly
controlled $K$-theory defined by Anderson-Munkholm in \cite{AM90}.
Finally, there are Atiyah-Hirzebruch spectral sequences computing
both the groups $\tK _{*+1}^{bc}(\xi ^+)$ (see \cite[Theorem
4.1]{AM90}) and $H_*^{\g}(T; \kz)$ (see \cite[Section 8]{Qu82}).  It
is easy to verify that the two spectral sequences are {\it
canonically identical}: they have the same $E^2$-terms and the same
differentials. Combining these results, and keeping track of the
various maps appearing in the sequence of isomorphisms, one can give
an alternate proof of the {\bf Fact}.

\vskip 5pt

\noindent {\bf (2)} We also point out that, from the $\g$-action on
the Bass-Serre tree $T$, it is easy to obtain constraints on the isomorphism
type of groups inside the collection $\V$.  Indeed, any such group must
be the stabilizer of a bi-infinite geodesic $\gamma \subset T$ (see {\bf Claim 2}),
and must act cocompactly on $\gamma$.  Recall that infinite
virtually cyclic subgroups are of two types: those that surject onto the
infinite dihedral group $D_\infty$, and those who don't.  

The groups that surject 
onto $D_\infty$ always split as an amalgamation $A*_CB$, with all three groups
$A,B,C$ finite, and $C$ of index two in both $A,B$.  Observe that if $V\in \V$ is
of this type, then under the action of $V$ on $\gamma\subset T$, $A$ and $B$ can be
identified with the stabilizers of a pair of vertices $v,w$, and $C$ can be 
identified with the stabilizer of the segment joining $v$ to $w$.  In particular,
$C$ {\it must be a subgroup of an edge stabilizer}, and hence is conjugate
(within $\g$) to a finite subgroup of $H$.  Furthermore, since $A,B$ both 
stabilize a pair of vertices, they must be conjugate (within $\g$) to a finite
subgroup of either $G_1$ or $G_2$. 

The groups that do not surject onto $D_\infty$ are of the form 
$F\rtimes_\alpha \mZ$, where $F$ is a finite group, and $\alpha \in Aut(F)$
is an automorphism.  If $V\in \V$ is of this form, then for the action of $V$
on $\gamma \subset T$, one has that $F$ can be identified with the subset
of $V$ that pointwise fixes $\gamma$, while the $\mZ$ component acts
on $\gamma$ via a translation.  In particular, $F$ is again conjugate (within
$\g$) to a finite subgroup of $H$.

In particular, if we are given an explicit amalgamation $\g= G_1*_HG_2$,
and we have knowledge of the finite subgroups inside the groups $H,G_1,G_2$,
then we can readily identify up to isomorphism the possible groups arising
in the collection $\V$.  If one has knowledge of the Nil-groups associated
with these various groups, our Main Theorem can be used to get corresponding
information about the Waldhausen Nil-group associated to $\g$.

\section{Concluding Remarks}

Having completed the proof of our Main Theorem, we now proceed to
isolate a few interesting corollaries.  As mentioned earlier, from
the viewpoint of topological applications, the most interesting
situation is the case where $R=\mZ$, i.e. integral group rings.

\begin{Cor}
Let $\g= G_1*_HG_2$ be an amalgamation, and assume that the
Farrell-Jones Isomorphism Conjecture holds for the groups $\g,
G_1,G_2$, and that $H$ is finite.  Then the associated
Waldhausen Nil-group $Nil^W_*(\mZ H; \mZ[G_1-H], \mZ[G_2-H])$ is
either trivial, or an infinitely generated torsion group.
\end{Cor}

\begin{Prf}
Note that since $H$ is finite, the amalgamation is acylindrical, and
our Main Theorem applies.  So the Waldhausen group we are interested in
splits as a direct sum of Nil-groups associated to a
particular collection $\V$ of virtually cyclic subgroups.  It is well known
that the Nil-groups associated to virtually cyclic groups are either trivial
or infinitely generated (see \cite{F77}, \cite{G1}, \cite{R}).  Furthermore, 
these groups are known to be purely torsion (see \cite{We81}, \cite{CP02}, 
\cite{KT03}, \cite{G2}), giving us the second statement.
\end{Prf}

A special case of the above corollary is worth mentioning:

\begin{Cor}
Let $\g= G_1*_HG_2$ be an amalgamation, and assume that 
$G_1,G_2$, and $H$ are all finite.  Then the associated
Waldhausen Nil-group $Nil^W_*(\mZ H; \mZ[G_1-H], \mZ[G_2-H])$ is
either trivial, or an infinitely generated torsion group.
\end{Cor}

\begin{Prf}
$G_1,G_2$ trivially satisfy 
the Farrell-Jones isomorphism conjecture, as they are finite.  Furthermore, the group
$\g$ is $\delta$-hyperbolic, and hence by a recent result of
Bartels-L\"uck-Reich \cite{BLR}, also satisfies the isomorphism
conjecture.  Hence the hypotheses of our previous corollary are 
satisfied.
\end{Prf}

Finally, we observe that the Bartels-L\"uck-Reich result \cite{BLR}
establishes the Isomorphism conjecture in all dimensions $i$ and
arbitrary rings with unity $R$.  In particular, the hypotheses of
our Main Theorem hold for arbitrary amalgamations of finite groups, giving:

\begin{Cor}
Let $\g=G_1*_HG_2$ be an amalgamation, where $G_1,G_2$, and $H$ are
all finite.  Then, for arbitrary rings with unity $R$, we have isomorphisms:
$$Nil^W_*(R H; R[G_1-H], R[G_2-H])\cong \bigoplus _{V\in \V}
H_*^{V}(E_{\fin}V\rightarrow *; \kr)$$ where
$H_*^{V}(E_{\fin}V\rightarrow *;\kr)$ are the cokernels of the
relative assembly maps associated to the virtually cyclic subgroups
$V\in \V$, and the collection $\V$ is as in the statement of our Main Theorem.
\end{Cor}

We point out that the special case of the modular group $\g=PSL_2(\mZ)=\mZ_2*\mZ_3$
has also been independently studied by Davis-Khan-Ranicki (see \cite[Section 3.3]{DKR}).

\end{document}